\documentclass[titlepage,12pt]{article} 
\usepackage[title]{appendix}
\usepackage{amssymb,amsthm,amsmath} 
\usepackage[a4paper]{geometry}
\usepackage{datetime2}
\usepackage[utf8]{inputenc}
\usepackage{hyperref}

\date{}

\newcommand{\z}{\mathbb{Z}}

\newcommand{\re}{\mathbb{R}}

\newcommand{\ep}{\varepsilon}

\newcommand{\dint}{\int\!\!\!\!\int}
\newcommand{\ud}{u_{\delta}}
\newcommand{\HH}{\mathcal{H}}
\newcommand{\ip}{\Lambda_{0,p}}
\newcommand{\dzp}{\delta\to 0^{+}}
\newcommand{\glim}{\Gamma\mbox{--}\lim}
\newcommand{\gliminf}{\Gamma\mbox{--}\liminf}
\newcommand{\ld}{\Lambda_{\delta,p}}
\newcommand{\ldp}{\Lambda_{\delta,p}}
\newcommand{\ldph}{\widehat{\Lambda}_{\delta,p}}
\newcommand{\lzp}{\Lambda_{0,p}}
\newcommand{\sd}{S_{\delta}}
\newcommand{\sn}{\mathbb{S}^{d-1}}
\newcommand{\pca}{\mathcal{PC\kern-1pt A}}

\renewcommand{\subparagraph}[1]{\paragraph{\textmd{\textit{#1}}}
}

\newtheorem{thm}{Theorem}[section]

\title{On the shape factor of interaction laws for a non-local approximation of the Sobolev norm and the total variation}

\author{Clara Antonucci\vspace{1ex}\\ 
{\normalsize Scuola Normale Superiore} \\
{\normalsize PISA (Italy)}\\
{\normalsize e-mail: \texttt{clara.antonucci@sns.it}}
\and
Massimo Gobbino\vspace{1ex}\\ 
{\normalsize Universit\`a degli Studi di Pisa} \\
{\normalsize PISA (Italy)}\\  
{\normalsize e-mail: \texttt{massimo.gobbino@unipi.it}}
\and
Matteo Migliorini\vspace{1ex}\\ 
{\normalsize Scuola Normale Superiore} \\
{\normalsize PISA (Italy)}\\
{\normalsize e-mail: \texttt{matteo.migliorini@sns.it}}
\and
Nicola Picenni\vspace{1ex}\\ 
{\normalsize Scuola Normale Superiore} \\
{\normalsize PISA (Italy)}\\
{\normalsize e-mail: \texttt{nicola.picenni@sns.it}}
}

\begin{document}
\maketitle
\begin{abstract}

We consider the family of non-local and non-convex functionals introduced in a recent paper by H.~Brezis and H.-M.~Nguyen. These functionals Gamma-converge to a multiple of the Sobolev norm or the total variation, depending on a summability exponent, but the exact values of the constants are unknown in many cases.

We describe a new approach to the Gamma-convergence result that leads in some special cases to the exact value of the constants, and to the existence of smooth recovery families.
	
\vspace{6ex}

\noindent{\bf Mathematics Subject Classification 2010 (MSC2010):}
26B30, 46E35.

% 26B30   	Absolutely continuous functions, functions of bounded variation
% 46E35   	Sobolev spaces and other spaces of "smooth'' functions, embedding theorems, trace theorems
\vspace{6ex}


\vspace{6ex}

% \listoftodos

\end{abstract}

%%%%%%%%%%%%%%%%%%%%%
%                   %
%   Inizio lavoro   %
%                   %
%%%%%%%%%%%%%%%%%%%%%

\section{Introduction}

In the recent paper~\cite{2018-AnPDE-BreNgu}, H.~Brezis and H.-M.~Nguyen introduced the family of non-local functionals
\begin{equation}
\ldp(\varphi,u,\Omega):=\dint_{\Omega^{2}}\varphi\left(\frac{|u(y)-u(x)|}{\delta}\right)\frac{\delta^{p}}{|y-x|^{d+p}}\,dx\,dy,
\label{defn:idphi}
\end{equation}
where $d$ is a positive integer, $\Omega\subseteq\re^{d}$ is an open set, $\delta>0$ is a real parameter, $p\geq 1$ is a real number, $u:\Omega\to\re$ is a measurable function, and $\varphi:[0,+\infty)\to[0,+\infty)$ is a measurable function that describes the extent to which a pair $(x,y)\in\Omega^{2}$ contributes to the double integral (\ref{defn:idphi}). For this reason, in the sequel we call $\varphi$ the ``interaction law''. A general enough class of admissible interaction laws is the set $\mathcal{A}$ of all functions $\varphi:[0,+\infty)\to[0,+\infty)$ that are not identically equal to zero and such that
\begin{itemize}

\item  $\varphi$ is nondecreasing and lower semicontinuous on $[0,+\infty)$, and actually continuous except at a finite number of points in $(0,+\infty)$,

\item  there exists a constant $a$ such that $\varphi(t)\leq at^{2}$ for every $t\in[0,1]$,

\item  there exists a constant $b$ such that $\varphi(t)\leq b$ for every $t\geq 0$.

\end{itemize}

The basic example is when $\varphi(x)$ coincides with
\begin{equation}
\renewcommand{\arraystretch}{1}
\varphi_{k}(x):=\left\{
\begin{array}{l@{\qquad}l}
 0 & \mbox{if } x\in[0,k],  \\
 1 & \mbox{if } x>k,   
\end{array}
\right.
\label{defn:phi-1}
\end{equation}
where $k$ is a positive real number.

The asymptotic behavior of the family $\ldp$, starting with the model case where $\varphi=\varphi_{1}$, was investigated in a series of papers~\cite{2006-CRAS-BouNgu,2015-Lincei-Brezis,2017-CRAS-BreNgu,2018-AnPDE-BreNgu,2006-JFA-Nguyen,2007-CRAS-Nguyen,2008-JEMS-Nguyen,2011-Duke-Nguyen}. The general idea is that $\ldp(\varphi,u,\Omega)$ is proportional, in the limit as $\delta\to 0^{+}$, to the functional
\begin{equation}
\ip(u,\Omega):=\left\{
\begin{array}{l@{\qquad}l}
    \displaystyle\int_{\Omega}|\nabla u(x)|^{p}\,dx  & \mbox{if $p>1$ and $u\in W^{1,p}(\Omega)$,}   \\
    \mbox{total variation of $u$ in $\Omega$}  & \mbox{if $p=1$ and $u\in BV(\Omega)$,}  \\    +\infty & \mbox{otherwise.}
\end{array}\right.
\end{equation}

The first result in this direction concerns the pointwise limit, at least in the case of smooth functions with compact support. Indeed, for every admissible interaction law $\varphi\in\mathcal{A}$ it turns out that
\begin{equation}
\lim_{\dzp}\ld\left(\varphi,u,\re^d\right)=G_{d,p}\cdot N(\varphi)\cdot\ip\left(u,\re^{d}\right)
\qquad\quad
\forall u\in C^{1}_{c}(\re^{d}),
\label{thmbibl:pointwise}
\end{equation}
where $G_{d,p}$ is a geometric constant, and $N(\varphi)$ is a normalization constant, which we call the ``scale factor'' of the interaction law $\varphi$. These two constants are defined as ($v$ is any element of the unit sphere $\sn$ in $\re^{d}$)
\begin{equation}
G_{d,p}:=\frac{1}{p}\int_{\sn}|\langle v,\sigma\rangle|^{p}\,d\sigma,
\hspace{4em}
N(\varphi):=\int_{0}^{+\infty}\frac{\varphi(t)}{t^{2}}\,dt.
\end{equation}

The equality in (\ref{thmbibl:pointwise}) holds true also for every $u\in W^{1,p}(\re^{d})$ if $p>1$ (but not necessarily if $p=1$).

The surprise comes with the Gamma-limit, which is expected to be of the form
\begin{equation}
\glim_{\dzp}\ld\left(\varphi,u,\re^d\right)=G_{d,p}\cdot N(\varphi)\cdot K_{d,p}(\varphi)\cdot\ip\left(u,\re^{d}\right)
\qquad
\forall u\in L^{p}(\re^{d}),
\label{thmbibl:g-conv}
\end{equation}
where $K_{d,p}(\varphi)\in(0,1]$ is a suitable constant, whose appearance was defined in~\cite{2018-AnPDE-BreNgu} ``mysterious and somewhat counterintuitive''. This result was proved in~\cite{2011-Duke-Nguyen} in the special case $\varphi=\varphi_{1}$ with general exponent $p\geq 1$, and in~\cite{2018-AnPDE-BreNgu} for general $\varphi\in\mathcal{A}$ but special exponent $p=1$. As far as we know, the case with general interaction law $\varphi\in A$ and general exponent $p\geq 1$ has never been written explicitly, even if a paper in this direction was announced in~\cite{2018-AnPDE-BreNgu}.

The constant $K_{d,p}(\varphi)$ is invariant by both horizontal and vertical rescaling, namely it does not change when we replace $\varphi(t)$ with $\alpha\varphi(\beta t)$ for some positive constants $\alpha$ and $\beta$. For this reason we call it the ``shape factor'' of the interaction law $\varphi$. 

Computing shape factors is a difficult task, even in dimension one (and actually we think that they never depend on $d$). In this note we describe a new approach to the Gamma-convergence problem that was carried out in our recent papers~\cite{AGMP:log-2,AGP:log-e}, and allowed to compute the shape factor of some special interaction laws. In the sequel
\begin{itemize}

\item $\mathcal{A}_{0}$ denotes the set of all $\varphi\in\mathcal{A}$ such that $\varphi(t)=0$ for every $t\in[0,1]$,

\item $\pca$ denotes the set of interaction laws that can be written in the form
\begin{equation}
\varphi(t)=\sum_{k=1}^{m}\lambda_{k}\varphi_{k}(t)
\qquad
\forall t\geq 0
\label{defn:pca}
\end{equation}
for some positive integer $m$, and some nonnegative coefficients $\lambda_{1},\ldots,\lambda_{m}$ (not all equal to 0),

\item  $\pca_{2}$ denotes the set of interaction laws of the form (\ref{defn:pca}) with coefficients equal in packages of powers of two, namely $\lambda_{2}=\lambda_{3}$, $\lambda_{4}=\ldots=\lambda_{7}$, $\lambda_{8}=\ldots=\lambda_{15}$, and so on.

\end{itemize}

Our first result (see Theorem~1.1 in~\cite{AGMP:log-2}) is the computation of the shape factor of the interaction laws defined in (\ref{defn:phi-1}). This settles a conjecture stated in~\cite{2007-CRAS-Nguyen,2011-Duke-Nguyen}.

\begin{thm}\label{thm:log-2}
Let us consider the interaction laws $\varphi_{k}(t)$ defined by (\ref{defn:phi-1}), with $k$ any positive real number. 

Then the Gamma-convergence result (\ref{thmbibl:g-conv}) holds true with
\begin{equation}
K_{d,p}(\varphi_{k}):=\left\{
\begin{array}{l@{\qquad}l}
    \displaystyle\frac{1}{p-1}\left(1-\frac{1}{2^{p-1}}\right)  & \mbox{if $p>1$,}   \\
    \log 2  & \mbox{if $p=1$.}  
\end{array}\right.
\end{equation}

\end{thm}

Two open questions raised in~\cite{2018-AnPDE-BreNgu} concerned the possibility that $K_{d,1}(\varphi)<1$ for every $\varphi\in\mathcal{A}$, and that $K_{d,1}(\varphi)$ might depend only on the behavior of $\varphi$ in a neighborhood of the origin. Our second result (see Theorem~1.3 and Theorem~1.4 in~\cite{AGP:log-e}) provides a negative answer to both questions.

\begin{thm}\label{thm:log-e}
Let us consider the Gamma-convergence result (\ref{thmbibl:g-conv}) in the case $p=1$.

Then in any space dimension $d$ it turns out that
\begin{equation}
\sup\{K_{d,1}(\varphi):\varphi\in\pca_{2}\}=\max\{K_{d,1}(\varphi):\varphi\in\mathcal{A}_{0}\}=1.
\end{equation}

\end{thm}

In the following sections we present the main steps of our strategy, without technicalities.  A notion that plays a central role in the sequel is what we call vertical $\delta$-segmentation. The vertical $\delta$-segmentation of a real-valued function $w(x)$ is the function $\sd w(x):=\delta\lfloor\delta^{-1}w(x)\rfloor$ that takes its values in $\delta\z$, and is uniquely characterized by the fact that $S_{\delta}w(x)=k\delta$ for some $k\in\z$ if and only if $k\delta\leq w(x)<(k+1)\delta$.

\setcounter{equation}{0}
\section{Aggregation/segregation problems}\label{sec:dino}

\subparagraph{Discrete setting}

Let us consider a positive integer $n$, a positive integer $k$, and a nonincreasing function $h:\{0,1,\ldots,n-1\}\to\re$, called ``hostility function''. A discrete arrangement is any function $u:\{1,\ldots,n\}\to\z$. For any such function, we define the total $k$-hostility as 
\begin{equation}
\HH_{k}(u):=\sum_{i=1}^{n}\sum_{j=i}^{n}\varphi_{k}(|u(j)-u(i)|)\cdot h(|j-i|),
\end{equation}
where $\varphi_{k}$ is the interaction law defined by (\ref{defn:phi-1}). Just to help intuition, we think of $u$ as an arrangement of $n$ dinosaurs placed in the points $\{1,\ldots,n\}$. There are different species of dinosaurs, indexed by integer numbers, so that $u(i)$ denotes the species of the dinosaur in position $i$. Two dinosaurs placed in the points $i$ and $j$ are hostile to each other if and only if the integers $u(i)$ and $u(j)$ representing their species differ by at least $k+1$, and in this case the real number $h(|j-i|)$ measures the ``hostility'' between the two dinosaurs. As expected, the closer are the positions of the dinosaurs, the larger is their hostility. 

The result is that the total $k$-hostility is minimized by monotone arrangements, namely those in which all dinosaurs of the same species are close to each other, and the groups corresponding to different species are sorted in ascending or descending order. More formally, for every discrete arrangement $u$, it turns out that $\HH_{k}(u)\geq\HH_{k}(Mu)$, where $Mu$ is the nondecreasing rearrangement of $u$, defined as the unique nondecreasing arrangement whose level sets have the same number of elements of the level sets of $u$. 

The proof of this quite intuitive statement seems to be somewhat non-trivial (see~\cite[Section~2.3]{AGMP:log-2}).

\subparagraph{Semi-discrete setting}

Let us consider now an interval $(a,b)\subseteq\re$, a positive integer $k$, and  a nonincreasing  hostility function $c:(0,b-a)\to\re$. A semi-discrete arrangement is any measurable function $u:(a,b)\to\z$ with finite image. For any such function $u$, we define the total $k$-hostility as
\begin{equation}
\mathcal{F}_{k}(c,u):=\dint_{(a,b)^{2}}\varphi_{k}(|u(y)-u(x)|)\cdot c(|y-x|)\,dx\,dy.
\label{defn:th-sd}
\end{equation}

Again the result is that a nondecreasing rearrangement does not increase the total $k$-hostility. The proof follows from the discrete counterpart through an approximation argument (any semi-discrete arrangement can be approximated with one whose level sets are finite unions of intervals with the same length). We refer to~\cite[Section~2.4]{AGMP:log-2} for the details. These results could also be deduced from (and are actually equivalent to) some rearrangement inequalities proved independently in the 1970s (see~\cite{1973-JCT-Taylor,1974-AnIF-GarRod}).

\setcounter{equation}{0}
\section{Estimating the Gamma-liminf from below}\label{sec:liminf}

\subparagraph{From dimension one to any space dimension}

Let us assume that there exists a constant $\Gamma_{0}$ such that
\begin{equation}
\gliminf_{\dzp}\ld\left(\varphi,u,\re^{d}\right)\geq G_{d,p}\cdot \Gamma_{0}\cdot\lzp\left(u,\re^{d}\right)
\qquad
\forall u\in L^{p}(\re^{d})
\label{hp:1-dim}
\end{equation}
when $d=1$. Then the same inequality holds true in any space dimension $d$.

The proof of this implication relies on an integral-geometric representation. The basic idea is that $\lzp(u,\re^{d})$ is the average of the analogous functional computed over all one-dimensional sections of $u$, namely over all restrictions of $u$ to lines. A similar representation holds true for $\ldp(\varphi,u,\re^{d})$. Moreover, the convergence $\ud\to u$ in $L^{p}(\re^{d})$ implies the convergence in $L^{p}(\re)$ of almost all one-dimensional sections. At this point, passing from dimension one to dimension~$d$ is just an application of Fatou's Lemma. 

We refer to~\cite[Section~4]{AGMP:log-2} for the details. This step requires no special assumption on $\varphi$.

\subparagraph{Localization technique}

Let us assume that there exists a constant $\Gamma_{0}$ such that, for every interval $(a,b)\subseteq\re$, and every family $\{\ud\}_{\delta>0}\subseteq L^{p}((a,b))$, it happens that
\begin{equation}
\liminf_{\dzp}\ld(\varphi,\ud,(a,b))\geq G_{1,p}\cdot \Gamma_{0}\cdot\frac{1}{(b-a)^{p-1}}\cdot\left(\liminf_{\dzp}\operatorname{osc}(\ud,(a,b))\right)^{p},
\label{hp:loc2glob}
\end{equation}
where $\operatorname{osc}(\ud,(a,b))$ denotes the essential oscillation of $\ud$ in $(a,b)$, namely the difference between the essential supremum and the essential infimum of $\ud$. Then (\ref{hp:1-dim}) holds true with $d=1$, and hence with any $d$.

The proof of this implication is quite classical. Given a family $\ud\to u$ in $L^{p}(\re)$, we approximate $u$ with any piecewise affine function $v$ whose graph is obtained by connecting points of the graph of $u$ corresponding to Lebesgue points of $u$. From estimate (\ref{hp:loc2glob}) applied in each interval, we deduce that the liminf of $\ld(\varphi,\ud,\re)$ is greater than or equal to $G_{1,p}\cdot \Gamma_{0}\cdot\lzp\left(v,\re\right)$, and we conclude by observing that $\lzp\left(u,\re\right)$ is the supremum of $\lzp\left(v,\re\right)$ as $v$ varies over all piecewise affine approximations of $u$. 

We refer to~\cite[Section~3.2]{AGMP:log-2} for the details. Also this step requires no special assumption on $\varphi$.

\subparagraph{Reduction to nondecreasing step functions}

Establishing (\ref{hp:loc2glob}) is the core of the proof of any estimate from below for the Gamma-liminf. It is also the point where assuming that $\varphi\in\pca$ yields a great simplification, because in this case the functional $\ldp$ is nonincreasing by vertical $\delta$-segmentation and monotone rearrangement. More precisely, let $\ud:(a,b)\to\re$ be any measurable function with bounded oscillation, let $\sd\ud$ be its vertical $\delta$-segmentation, and let $M\sd\ud$ be the nondecreasing rearrangement of $\sd\ud$. Then it turns out that 
\begin{equation}
\ldp(\varphi,\ud,(a,b))\geq\ldp(\varphi,M\sd\ud,(a,b)),
\label{ineq:MS}
\end{equation}
while the oscillation of $\ud$ is more or less the same as the oscillation of $M\sd\ud$ (their difference is less than $\delta$).  In proving (\ref{ineq:MS}) for $\varphi\in\pca$, due to the linearity of $\ldp$ with respect to $\varphi$, we can limit ourselves to the special case where $\varphi=\varphi_{k}$, in which case the functional $\ldp$ is equivalent to the semi-discrete total hostility $\mathcal{F}_{k}$ defined in~(\ref{defn:th-sd}) with hostility function $c(\sigma):=\delta^{p}\sigma^{-1-p}$, and hence the result follows from the corresponding estimate for the semi-discrete aggregation problem. We also show that
\begin{equation}
\liminf_{\dzp}\ldp(\varphi,\ud,(a,b))=\liminf_{\dzp}\ldph(\varphi,\ud,(a,b)),
\label{ldp-ldph}
\end{equation}
where $\ldph(\varphi,u,(a,b))$ is defined as $\ldp(\varphi,u,(a,b))$, just with the double integration over $(a,b)\times\re$ instead of $(a,b)\times(a,b)$, which leads to integrals that are simpler to compute explicitly. We refer to~\cite[Section~3.1]{AGMP:log-2} and~\cite[Lemma~4.1 and Lemma~4.2]{AGP:log-e} for the details.

\subparagraph{Reduction to the asymptotic study of multi-variable minimum problems}

In the previous paragraphs we have reduced ourselves to showing (\ref{hp:loc2glob}) for families of nondecreasing step functions $\ud:(a,b)\to\re$ with finite image contained in $\delta\z$. Up to vertical translations, any such function depends only on the lengths $\ell_{1},\ldots,\ell_{n}$ of the steps, where $n\sim\delta^{-1}\cdot\operatorname{osc}(\ud,(a,b))$ denotes the number of the steps. With a homothety we can also rescale $(a,b)$ to an interval of unit length, and by homogeneity we obtain that
\begin{equation}
\ldph(\varphi,\ud,(a,b))=\frac{\delta^{p}}{(b-a)^{p-1}}P_{n,\varphi,p}(\ell_{1},\ldots,\ell_{n}),
\nonumber
\end{equation}
for a suitable multi-variable function $P_{n,\varphi,p}$ defined as follows. We consider the representation of $\varphi$ in the form (\ref{defn:pca}), and then we set
\begin{equation}
P_{n,\varphi,p}:=\sum_{k=1}^{m}\lambda_{k}\sum_{i=1}^{n-k}\left(\int_{S_{i,k}}^{S_{i,k+1}}\frac{1}{\sigma^{p}}\,d\sigma+\int_{S_{i+1,k}}^{S_{i,k+1}}\frac{1}{\sigma^{p}}\,d\sigma\right),
\end{equation}
where $S_{i,h}:=\ell_{i}+\ldots+\ell_{i+h-1}$.
Setting
\begin{equation}
I_{n,p}(\varphi):=\inf\left\{P_{n,\varphi,p}(\ell_{1},\ldots,\ell_{n}):(\ell_{1},\ldots,\ell_{n})\in(0,+\infty)^{n},\ \ell_{1}+\ldots+\ell_{n}=1\right\},
\label{defn:inf}
\end{equation}
we can prove that
\begin{equation}
\liminf_{\dzp}\ldph(\varphi,\ud,(a,b))\geq\frac{1}{(b-a)^{p-1}}\cdot\liminf_{\dzp}(\delta n)^{p}\cdot\liminf_{n\to +\infty}\frac{I_{n,p}(\varphi)}{n^{p}}.
\label{3-liminf}
\end{equation}

Recalling that $\delta n\sim\operatorname{osc}(\ud,(a,b))$, the proof of an estimate of the form (\ref{hp:loc2glob}) has been reduced to the asymptotic study of a family of multi-variable minimum problems.

\subparagraph{Special piecewise constant interaction laws}
 
The minimum problems (\ref{defn:inf}) can be very complicated, but in some special cases they are quite simple. For example, when $\varphi=\varphi_{1}$ and $p=1$ it turns out that 
\begin{equation}
P_{n,\varphi_{1},1}(\ell_{1},\ldots,\ell_{n}) = \log\frac{(\ell_{1}+\ell_{2})^{2}}{\ell_{1}\ell_{2}}+\log\frac{(\ell_{2}+\ell_{3})^{2}}{\ell_{2}\ell_{3}}+\ldots+\log\frac{(\ell_{n-1}+\ell_{n})^{2}}{\ell_{n-1}\ell_{n}}.
\nonumber
\end{equation}

All the fractions are greater than or equal to 4, and hence $I_{n,1}(\varphi_{1})=(n-1)\log 4$, with the minimum realized when all the variables are equal. Since $G_{1,1}=2$, we obtain that (\ref{hp:loc2glob}) holds true with $\Gamma_{0}=\log 2$, and this leads to the proof of Theorem~\ref{thm:log-2} in the case $p=1$. An analogous inequality settles also the case $p>1$ (see~\cite[Proposition~3.2]{AGMP:log-2}).

More generally, we can handle interaction laws $\varphi\in\pca_{2}$. Indeed, in this case a telescopic effect (see~\cite[Section~3]{AGP:log-e}) simplifies the computation of $P_{n,\varphi,1}$, and we obtain that $I_{n,1}(\varphi)\sim n\cdot 2\log 2\cdot(\lambda_{1}+\lambda_{2}+\lambda_{4}+\ldots+\lambda_{2^{m-1}})$. 

In the case where all coefficients are equal, letting $m\to +\infty$ we find a sequence of piecewise constant interaction laws whose shape factors tend to~1, which proves that the supremum of $K_{d,1}(\varphi)$ in $\pca_{2}$ is one, namely the first conclusion of Theorem~\ref{thm:log-e}. Then we consider the interaction law $\theta(t):=\min\{1,\max\{t-1,0\}\}$, and we approximate it from below by a sequence of rescalings of interaction laws in $\pca_{2}$ with shape factors that tend to one. This proves the second conclusion of Theorem~\ref{thm:log-e}. We refer to \cite[Section~5]{AGP:log-e} for the details.

\setcounter{equation}{0}
\section{Estimating the Gamma-limsup from above}

\subparagraph{Reduction to subsets that are dense in energy}

It is well-known that it is enough to show the existence of recovery families for every $u$ belonging to a subset that is dense in energy with respect to the limit functional. Good classes in this case are smooth or piecewise affine functions with compact support. 

\subparagraph{Vertical $\delta$-segmentation provides a recovery family in dimension one}

Let us consider the special case $\varphi=\varphi_{1}$. Then for every piecewise $C^{1}$ function $u$ with compact support, the family $\{\sd u\}$ of its vertical $\delta$-segmentations is a recovery family for $u$. This is proved in~\cite[Section~3.3]{AGMP:log-2}, and the argument (based on the dominated convergence theorem for integrals) is similar to the proof of the pointwise convergence result (\ref{thmbibl:pointwise}). In~\cite{AGP:log-e} it is explained how the result can be extended to larger classes of interaction laws.
 
Analogous arguments should lead to a proof that $\{\sd u\}$ is a recovery family for any piecewise smooth function $u$ whenever the interaction law $\varphi$ is in $\pca$ and the infimum in (\ref{defn:inf}) is realized (at least asymptotically) when $\ell_{1}=\ldots=\ell_{n}$.

\subparagraph{Smooth recovery families in dimension one}

Let us assume that, for some $\varphi\in\mathcal{A}$, the vertical $\delta$-segmentation provides a recovery family for all piecewise affine functions with compact support in dimension one. Then for every $u\in L^{p}(\re)$ there exists a recovery family in $C^{\infty}_{c}(\re)$.

It is enough to prove this result when $u$ is piecewise affine with compact support. In this case we already know that $\{\sd u\}$ is a recovery family for $u$. For every fixed $\delta>0$, the function $\sd u$ has a very simple structure, namely it is a step function with a finite number of steps, and level sets are a finite number of intervals. For any such function, there exists a family $\{u_{\ep,\delta}\}_{\ep>0}$ of functions in $C^{\infty}_{c}(\re)$ such that
\begin{equation}
\lim_{\ep\to 0^{+}}u_{\ep,\delta}=\sd u \quad\mbox{in }L^{p}(\re)
\qquad\mbox{and}\qquad
\lim_{\ep\to 0^{+}}\ldp(\varphi,u_{\ep,\delta},\re)=\ldp(\varphi,\sd u,\re).
\nonumber
\end{equation}

To this end, it is enough to round the corners off by defining $u_{\ep,\delta}$ as the convolution of $\sd u$ with a rescaled smooth kernel with compact support. We refer to~\cite[Section~3.4]{AGMP:log-2} for the details.

\subparagraph{From dimension one to any dimension}

Let us assume that, for some $\varphi\in\mathcal{A}$, the vertical $\delta$-segmentation provides a recovery family in dimension one. Then the same is true in any space dimension.

Indeed, the key observation is that the vertical $\delta$-segmentation commutes with taking one-dimensional sections. At this point, the conclusion follows from Fatou's Lemma and from the integral-geometric representation described at the beginning of Section~\ref{sec:liminf}. We refer to~\cite[Section~4]{AGMP:log-2} for the details.

\subparagraph{Smooth recovery families in any dimension}

Under the same assumptions of the previous step, smooth recovery families do exist in any space dimension. The basic idea is analogous to the one dimensional case. If $u$ is piecewise affine with compact support in $\re^{d}$, then $\{\sd u\}$ is a recovery family. For every fixed $\delta>0$, the function $\sd u$ is a step function whose level sets are the union of a finite number of polytopes. These functions can be approximated in energy by smooth functions (through convolution). 

\setcounter{equation}{0}
\section{Open problems}\label{sec:open}

\subparagraph{Semi-discrete aggregation problem in any space dimension}

Even if not related directly to the theory of the functionals considered in this paper, it could be interesting to state the semi-discrete aggregation problem for arrangements $u:\Omega\to\z$, where $\Omega$ is a bounded open set in $\re^{d}$. Finding the minimizers of the total hostility in this context seems to be a challenging problem in geometric measure theory, since it is not clear what plays the role of monotone arrangements in higher dimension.

\subparagraph{Gamma-limit for piecewise constant interaction laws}

We suspect that the second $\liminf$ in the right-hand side of (\ref{3-liminf}) is actually a limit. We also suspect that, for every $\varphi\in\pca$, the right-hand side of (\ref{3-liminf}) is actually the Gamma-limit of $\ldp(\varphi,u,(a,b))$, and not just an estimate from below for the Gamma-liminf. 

This would reduce the computation of shape factors of interaction laws in $\pca$ to the computation of the asymptotic behavior of the minimum problems (\ref{defn:inf}). 

\subparagraph{From piecewise constant to general interaction laws}

We suspect that the Gamma-limit of $\ld(\varphi,u,\re^{d})$ is the supremum of the Gamma-limits of $\ld(\psi,u,\re^{d})$, where $\psi\leq\varphi$ varies in the set of all piecewise constant interaction laws with steps of equal horizontal length. Any such $\psi$ is the rescaling of an element of $\pca$, and rescaling preserves the shape factor.
A confirmation of this conjecture would open the way for answering several questions raised in~\cite{2018-AnPDE-BreNgu}: a simplified proof of the Gamma-convergence result in full generality, a less implicit formula for shape factors, and existence of smooth recovery families.

%\bibliographystyle{MaxNew}
%\bibliography{../../../BibTeX/Nguyen}

\label{NumeroPagine}

\end{document}